\documentclass[11pt,a4paper]{article}
\usepackage[T1]{fontenc}
\usepackage[utf8]{inputenc}
\usepackage{lmodern,amsmath,amssymb,amsthm,mathtools}
\usepackage[margin=28mm]{geometry}
\usepackage{microtype}
\usepackage[hidelinks]{hyperref}
\hypersetup{pdftitle={Fefferman--Stein Type Distribution Function Inequalities on Product Spaces},pdfauthor={Jiecheng Chen}}
\numberwithin{equation}{section}
\newtheorem{theorem}{Theorem}[section]
\newtheorem{lemma}[theorem]{Lemma}
\newtheorem{proposition}[theorem]{Proposition}
\newtheorem{corollary}[theorem]{Corollary}
\theoremstyle{remark}
\newtheorem{remark}[theorem]{Remark}
\newcommand{\R}{\mathbb R}
\newcommand{\C}{\mathbb C}
\newcommand{\HH}{\mathcal H}
\newcommand{\Ed}{\mathcal E_d}
\newcommand{\Enm}{\mathcal E_{n,m}}
\newcommand{\ind}{\mathbf 1}
\newcommand{\Ph}{\Phi_\lambda}
\newcommand{\Th}{\Theta_\lambda}
\newcommand{\grad}{\nabla^{\perp}}
\DeclareMathOperator{\supp}{supp}
\allowdisplaybreaks[2]
\title{Fefferman--Stein Type Distribution Function Inequalities\\on Product Spaces}
\author{Jiecheng Chen\\\small Department of Mathematics, Zhejiang Normal University}
\date{}
\begin{document}
\maketitle
\begin{abstract}
We establish on $\R^n\times\R^m$ a Fefferman--Stein type distribution function inequality in which the nontangential maximal function of the product Poisson extension is controlled by the Lusin area integral function:
\[
 |\{N>\lambda\}|\le C\left(|\{A>\lambda\}|+\lambda^{-2}\int_{\{A\le\lambda\}}A^2\,dx\right).
\]
The central idea is to establish a local $L^2$ estimate for the radial maximal function. This yields a new proof of the corresponding one-parameter inequality, avoiding Green's formula and the treatment of sawtooth regions, and extends to the multiparameter setting.
\end{abstract}


\section{Introduction and main results}
Maximal functions and square functions are basic tools in the study of Hardy space theory. Comparing these functions is a fundamental problem in this theory, particularly in connection with the following two Fefferman--Stein type distribution function inequalities for the area integral function $A$ and the nontangential maximal function $N$:
\begin{align}
 |\{A>\lambda\}|&\le C|\{N>\lambda\}|+\frac{C}{\lambda^2}\int_{\{N\le\lambda\}}N^2\,dx,\label{eq:forward}\\
 |\{N>\lambda\}|&\le C|\{A>\lambda\}|+\frac{C}{\lambda^2}\int_{\{A\le\lambda\}}A^2\,dx.\label{eq:reverse}
\end{align}
Both inequalities have been studied systematically in the one-parameter setting. In the multiparameter setting, however, \eqref{eq:forward} has been extensively investigated, whereas progress on \eqref{eq:reverse} has remained elusive. In this paper we study the latter inequality on product spaces and establish \eqref{eq:reverse}.

In the one-parameter theory, distribution function comparisons in both directions between the area integral function and the nontangential maximal function have long been classical. Burkholder and Gundy~\cite{BG} studied distribution function inequalities for the area integral using probabilistic and martingale methods, connecting localized distribution control with global integral estimates. Fefferman and Stein~\cite{FS} established the corresponding real-variable framework in their theory of $H^p$ spaces of several variables. By using Green's formula and a subtle treatment of sawtooth regions, they obtained distribution estimates in both directions, namely the one-parameter versions of \eqref{eq:forward} and \eqref{eq:reverse}. Subsequent developments included weighted comparisons, geometric extensions to Lipschitz domains, and sharper ratio estimates. The weighted theory of Gundy and Wheeden~\cite{GW}, Dahlberg's results on Lipschitz domains~\cite{Dahlberg}, and the ratio estimates of Fefferman, Gundy, Silverstein, and Stein~\cite{FGSS} illustrate the role of distribution comparisons in weighted theory, complex boundary geometry, and stronger integral control, respectively. Murai and Uchiyama~\cite{MU} studied good-$\lambda$ estimates, while the work of Ba\~nuelos and Moore~\cite{BM1,BM2} further addressed sharp comparisons and the density of the area integral.

The difficulty in the multiparameter setting arises because the two groups of spatial variables have independent time variables, while the area integral function involves mixed derivatives in the two groups of space-time variables. The $H^p$ theory on the bidisc developed by Malliavin and Malliavin~\cite{MM} and by Gundy and Stein~\cite{GS} was pioneering. In~\cite{GS}, Gundy and Stein first established \eqref{eq:forward} and a probabilistic version of \eqref{eq:reverse} on the bidisc. Subsequently, Merryfield~\cite{Merryfield}, Sato~\cite{Sato}, and Chen~\cite{ChenThesis} (also see Chen~\cite{Chen1999} for its proof strategy in the one-parameter case) independently obtained distribution estimates of type \eqref{eq:forward}, controlling the area integral function by the nontangential maximal function on various product domains;. Chen and Wang~\cite{CW} extended \eqref{eq:forward} to products of complete Riemannian manifolds with nonnegative Ricci curvature, while Li~\cite{Li} extended it to product Shilov boundaries. For the reverse distribution inequality \eqref{eq:reverse}, however, there has been no further breakthrough.

Let $d_1=n$, $d_2=m$, where $n,m\ge1$, and set $X=\R^{d_1}\times\R^{d_2}$. For $i=1,2$, let $\Delta_i$ denote the Laplacian in the $i$th group of variables, and define
\[
 P_t^{(i)}=e^{-t\sqrt{-\Delta_i}},\qquad
 u(y_1,t_1;y_2,t_2)=P_{t_1}^{(1)}P_{t_2}^{(2)}f(y_1,y_2).
\]
Write $z_i=(y_i,t_i)$, $\grad_i=(\nabla_{y_i},\partial_{t_i})$, $v_{d_i}=|B_{\R^{d_i}}(0,1)|$, and
\begin{equation}
 \Gamma_\alpha(x)=\Gamma_{1,\alpha_1}(x_1)\times\Gamma_{2,\alpha_2}(x_2),\quad
 \Gamma_{i,a}(x_i)=\{(y_i,t_i):|y_i-x_i|<at_i\}.\label{eq:cones}
\end{equation}
Vectors and tensors between the two groups of variables are equipped with their natural Hilbert norms (which are simply Euclidean norms). For $\alpha=(\alpha_1,\alpha_2)$ and $\beta=(\beta_1,\beta_2)\in(0,\infty)^2$, define
\begin{align}
 &N^+f(x)=\sup_{t_1,t_2>0}|u(x_1,t_1;x_2,t_2)|,\label{eq:radial}\\
 &N_\alpha f(x)=\sup_{(z_1,z_2)\in\Gamma_\alpha(x)}|u(z_1,z_2)|,\label{eq:nt}\\
 &A_\beta f(x)=\left(\iint_{\Gamma_\beta(x)}|\grad_1\grad_2u(y,t)|^2
       \prod_{i=1}^2\frac{t_i\,dy_i\,dt_i}{v_{d_i}(\beta_i t_i)^{d_i}}\right)^{1/2}.\label{eq:area}
\end{align}
For a nonnegative function $F$, write $d_F(\lambda)=|\{x:F(x)>\lambda\}|$, and set
\begin{equation}
 \Th(s)=\min\{s^2,\lambda^2\},\qquad \Ph(F)=\int_X\Th(F(x))\,dx.\label{eq:phi}
\end{equation}
We have the following result.
\begin{theorem}\label{th:main}
For every $f\in C_c^\infty(\R^n\times\R^m)$ and $\lambda>0$, we have
\begin{equation}
 \Ph(N_\alpha f)\le C\Ph(A_\beta f),\label{eq:main}
\end{equation}
where $C$ depends only on the dimensions and on $\alpha,\beta$, and $\Ph$ is defined by \eqref{eq:phi}. The corresponding one-parameter result holds for Hilbert-valued functions, with a constant independent of the target Hilbert space.
\end{theorem}
Theorem~\ref{th:main} implies the reverse distribution estimate \eqref{eq:reverse} and the corresponding one-sided strong and weak $L^p$ estimates for $0<p<2$. The central idea of its proof is to establish the following local $L^2$ estimate for the radial maximal function, thereby reducing the study of the nontangential maximal function to that of the radial maximal function.

Let $M_i$ be the centered Hardy--Littlewood maximal operator on the $i$th spatial factor $\R^{d_i}$, and let $M_1M_2$ denote their composition in the indicated order.
\begin{theorem}\label{th:local}
There exists $\delta\in(0,1)$, depending only on $n,m,\beta$, such that for every measurable $E\subset X$ with $|E^c|<\infty$, the set
\[
 G_E=\{x:M_1M_2\ind_{E^c}(x)\le\delta\}
\]
satisfies
\begin{equation}
 |G_E^c|\le C|E^c|,\qquad
 \int_{G_E}(N^+f)^2\,dx\le C\int_E(A_\beta f)^2\,dx.\label{eq:local}
\end{equation}
This holds for every $f\in L^2(X)$. The corresponding one-parameter Hilbert-valued result also holds.
\end{theorem}

\section{One-parameter preliminaries}
This section presents the one-parameter results needed below.
\begin{lemma}\label{le:trunc}
Suppose $c_1,\ldots,c_k\ge0$. For all $s_1,\ldots,s_k\ge0$, we have
\begin{equation}
 \Th\left(\sum_{j=1}^kc_js_j\right)\le k\sum_{j=1}^k\max\{1,c_j^2\}\Th(s_j).\label{eq:trunc}
\end{equation}
If nonnegative functions $F_j\to F$ in $L^2(X)$, then $\Ph(F_j)\to\Ph(F)$ for every $\lambda>0$.
\end{lemma}
\begin{proof}
By the Cauchy--Schwarz inequality, $(\sum_jc_js_j)^2\le k\sum_j(c_js_j)^2$. Considering separately the case where some $c_js_j\ge\lambda$ and the case where all $c_js_j<\lambda$, we obtain
\[
 \min\left\{k\sum_j(c_js_j)^2,\lambda^2\right\}
 \le k\sum_j\min\{(c_js_j)^2,\lambda^2\}.
\]
The first assertion follows from $\Th(cs)\le\max\{1,c^2\}\Th(s)$. Also,
\begin{equation}
 |\Th(a)-\Th(b)|\le|a-b|(a+b),\qquad a,b\ge0.\label{eq:continuity}
\end{equation}
Integrating and applying Cauchy--Schwarz inequality gives the final assertion.
\end{proof}
For every dimension $d\ge1$, fix a radial function $\varphi_d\in C_c^\infty(B(0,1/4))$ with $\varphi_d\ge0$, $\int\varphi_d=1$, and $\varphi_d>0$ on $B(0,1/8)$. Set
\begin{equation}
 \psi_d(x)=\varphi_d(x)-2^{-d}\varphi_d(x/2).\label{eq:psi}
\end{equation}
Then $\int\psi_d=0$ and $\supp\psi_d\subset B(0,1/2)$. Choose also a nonnegative $\chi_d\in C_c^\infty(B(0,1))$ with $\int\chi_d=1$. Write $\psi_{d,t}(x)=t^{-d}\psi_d(x/t)$ and $\chi_{d,t}(x)=t^{-d}\chi_d(x/t)$.

Let $\HH$ be a Hilbert space (it may be a Euclidean space, whose dimension depends on the gradient and differs between the one- and two-parameter settings, or an $L^2$ space. For simplicity, we use the same notation throughout). Define
\begin{equation}
 \begin{split}
 &\Ed(\HH)=L^2\bigl(\R^d\times(0,\infty),t\,dx\,dt;\HH\bigr),\\
 &\|H\|_{\Ed(\HH)}=\left(\int_0^\infty\int_{\R^d}t\|H(x,t)\|_\HH^2\,dx\,dt\right)^{1/2}.
 \end{split}\label{eq:energy}
\end{equation}
Write $\psi_t=\psi_{d,t}$ and $\chi_t=\chi_{d,t}$. Define
\begin{equation}
 T_sH=\int_s^\infty\psi_{t-s}*H_t\,dt,\quad
 E_sH=T_sH-\chi_s*SH,\quad
 SH=\int_0^\infty\psi_t*H_t\,dt,\label{eq:operators}
\end{equation}
where $H_t(x)=H(x,t)$. The following maximal lemma for $T_s$ is crucial.
\begin{lemma}\label{le:max}
There exists a nonnegative $k\in L^1(\R^d)$ such that, for $H\in\Ed(\HH)$,
\begin{align}
 &\|SH\|_2^2\le C\|H\|_{\Ed(\HH)}^2,\label{eq:Sbound}\\
 &\sup_{s>0}\|E_sH(x)\|_\HH^2\le C\int_0^\infty t\bigl(K_t*\|H_t\|_\HH^2\bigr)(x)\,dt,\label{eq:Ebound}
\end{align}
where $K_t(x)=t^{-d}k(x/t)$. Consequently,
\begin{equation}
 \left\|\sup_{s>0}\|T_sH(\cdot)\|_\HH\right\|_{L^2(\R^d)}^2
 \le C\|H\|_{\Ed(\HH)}^2.\label{eq:Tbound}
\end{equation}
The constants are independent of the target Hilbert space $\HH$ and of $H$.
\end{lemma}
\begin{proof}
By a limiting argument, we may assume that $\supp H$ is compact in $\R^d\times(0,\infty)$. Since $\psi$ is radial, compactly supported, and has mean zero,
\[
 |\widehat\psi(\xi)|\le C\min\{|\xi|,|\xi|^{-1}\}.
\]
For $f\in L^2(\R^d;\HH)$, Plancherel's theorem gives
\[
 \int_0^\infty\|\psi_t*f\|_2^2\,\frac{dt}{t}
 =\int_{\R^d}\|\widehat f(\xi)\|_\HH^2
 \left[\int_0^\infty|\widehat\psi(t\xi)|^2\,\frac{dt}{t}\right]d\xi
 \le C\|f\|_2^2.
\]
By duality,
\[
 |\langle SH,f\rangle|
 \le\left[\int_0^\infty t\|H_t\|_2^2\,dt\right]^{1/2}
 \left[\int_0^\infty\|\widetilde\psi_t*f\|_2^2\,\frac{dt}{t}\right]^{1/2}
 \lesssim\|H\|_{\Ed(\HH)}\|f\|_2,
\]
where $\widetilde\psi(x)=\psi(-x)$. This proves \eqref{eq:Sbound}.

We next prove \eqref{eq:Ebound}. We can write
\begin{equation}
 E_sH(x)=\int_0^\infty(e_{s,t}*H_t)(x)\,dt,\label{eq:Ekernel}
\end{equation}
where $e_{s,t}=\ind_{\{t>s\}}\psi_{t-s}-\chi_s*\psi_t$. Choose a sufficiently large constant $C_0$ and put $a_r(x)=r^{-d}\ind_{\{|x|<C_0r\}}(x)$. Taylor's formula and the compact supports of $\psi,\chi$ yield
\begin{align*}
 &|e_{s,t}(x)|\le C(t/s)a_s(x),&&0<t<s/2,\\
 &|e_{s,t}(x)|\le C(s/t)a_t(x),&&t>2s,\\
 &|\chi_s*\psi_t(x)|\le Ca_s(x),&&s/2\le t\le2s,\\
 &|\psi_{t-s}(x)|\le Ca_{t-s}(x),&&s<t<2s.
\end{align*}
Accordingly, decompose $E_sH(x)$ into four terms:
\begin{align*}
 E_sH(x)&=\int_0^{s/2}(e_{s,t}*H_t)(x)\,dt
       +\int_{2s}^\infty(e_{s,t}*H_t)(x)\,dt\\
 &\quad-\int_{s/2}^{2s}(\chi_s*\psi_t*H_t)(x)\,dt
       +\int_s^{2s}(\psi_{t-s}*H_t)(x)\,dt\\
 &=I_1(x)+I_2(x)+I_3(x)+I_4(x).
\end{align*}
Set $F_t(x)=\|H_t(x)\|_\HH^2$. Cauchy--Schwarz inequality gives
\begin{align*}
 \|I_1(x)\|_\HH^2&\le C\int_0^{s/2}t\frac{t}{s}(a_s*F_t)(x)\,dt,\\
 \|I_2(x)\|_\HH^2&\le C\int_{2s}^\infty t\frac{s}{t}(a_t*F_t)(x)\,dt,\\
 \|I_3(x)\|_\HH^2&\le C\int_{s/2}^{2s}t(a_s*F_t)(x)\,dt.
\end{align*}
For $I_4$, we have
\begin{align*}
 \|I_4(x)\|_\HH^2
 &\le\left[\int_s^{2s}t\left(\frac{t-s}{t}\right)^\alpha(a_{t-s}*F_t)(x)\,dt\right]
 \left[\int_s^{2s}\left(\frac{t-s}{t}\right)^{-\alpha}\frac{dt}{t}\right]\\
 &\le C_\alpha\int_s^{2s}t\left(\frac{t-s}{t}\right)^\alpha(a_{t-s}*F_t)(x)\,dt,
\end{align*}
since
\[
 \int_s^{2s}\left(\frac{t-s}{t}\right)^{-\alpha}\frac{dt}{t}
 =\int_0^1u^{-\alpha}(1+u)^{\alpha-1}\,du=C_\alpha<\infty
\]
whenever $\alpha\in(0,1)$. Take
\[
 k(z)=C\left((1+|z|)^{-d-1}+|z|^{\alpha-d}\ind_{\{|z|<C_1\}}\right)\in L^1(\R^d).
\]
Estimate \eqref{eq:Ebound} follows from
\begin{align*}
 &\sup_{s>2t}\frac{t}{s}a_s(z)\le Ct\sup_{s>\max(2t,|z|/C_0)}s^{-d-1} \le Ct^{-d}(1+|z|/t)^{-d-1},\\
 &\sup_{0<s<t/2}\frac{s}{t}a_t(z)\le Ca_t(z)\le Ct^{-d}\ind_{\{|z|<C_0t\}},\\
 &\sup_{t/2\le s\le2t}a_s(z)\le Ct^{-d}\ind_{\{|z|<C_1t\}},\\
 &\sup_{0<r<t/2}\left(\frac{r}{t}\right)^\alpha a_r(z)\le Ct^{-\alpha}|z|^{\alpha-d}\ind_{\{|z|<C_1t\}}.
\end{align*}
Finally, to prove \eqref{eq:Tbound}, note that
\[
 \|\chi_s*SH(x)\|_\HH\le Cs^{-d}\int_{|x-y|<s}\|SH(y)\|_\HH\,dy
 \le CM(\|SH\|_\HH)(x),
\]
where $M$ is the Hardy--Littlewood maximal operator on $\R^d$. Combining this with \eqref{eq:Sbound}--\eqref{eq:Ebound} proves \eqref{eq:Tbound}.
\end{proof}
For every $\theta\in S^{d-1}$, radiality of $\psi$ and Fubini's theorem give
\begin{equation}
 \begin{split}
 c_d:=&\int_0^\infty e^{-r}\widehat\psi_d(r\theta)\,dr\\
=&\int_{\R^d}\varphi_d(y)\left(\frac{1}{1+(y\cdot\theta)^2}
       -\frac{1}{1+4(y\cdot\theta)^2}\right)dy>0.
 \end{split}\label{eq:cd}
\end{equation}
This value is independent of $\theta$. We have the following reconstruction formula.
\begin{lemma}\label{le:reconstruct}
Let $q\in L^2(\R^d;\HH)$ and $v(x,t)=P_tq(x)$, where $P_t$ is the Poisson semigroup on $\R^d$. For every $s>0$, the following identity holds in $L^2$:
\begin{equation}
 v(x,s)=-c_d^{-1}T_s(\partial_tv(\cdot,t))(x).\label{eq:reconstruct}
\end{equation}
\end{lemma}
\begin{proof}
For $\xi\ne0$, the Fourier multiplier on the right-hand side of \eqref{eq:reconstruct} is
\[
 -c_d^{-1}\int_0^\infty\widehat\psi_d(r\xi)(-|\xi|)e^{-(s+r)|\xi|}\,dr=e^{-s|\xi|}.
\]
Here we use the substitution $r|\xi|\mapsto r$ and the definition of $c_d$. Since $\partial_tv\in\Ed(\HH)$, Lemma~\ref{le:max} justifies \eqref{eq:reconstruct}.
\end{proof}
\begin{lemma}\label{le:aperture}
Let $F$ be a measurable function on $\R_+^{d+1}$ with values in a Hilbert space $\HH$, and, for $c>0$, define
\[
 S_cF(x)=\left(\int_{\Gamma_c(x)}\|F(y,t)\|_\HH^2\frac{t\,dy\,dt}{v_d(ct)^d}\right)^{1/2}.
\]
For fixed $a,b>0$, we have
\begin{equation}
 \Ph(S_aF)\le C\Ph(S_bF).\label{eq:aperture}
\end{equation}
\end{lemma}
\begin{proof}
The assertion is immediate if the right-hand side is infinite. Otherwise, set $E=\{S_bF\le\lambda\}$, $B=E^c$, and $G=\{M\ind_B\le\varepsilon\}$, where $M$ is the centered Hardy--Littlewood maximal operator. Its weak type $(1,1)$ bound gives $|G^c|\le C\varepsilon^{-1}|B|$. For any $(y,t)$ such that $G\cap B(y,at)\ne\varnothing$, choose $x\in G\cap B(y,at)$. Since $B(y,bt)\subset B(x,(a+b)t)$,
\[
 1-\frac{|E\cap B(y,bt)|}{|B(y,bt)|}\le C_{a,b}M\ind_B(x)\le C_{a,b}\varepsilon.
\]
Fix $\varepsilon$ sufficiently small that $|E\cap B(y,bt)|/|B(y,bt)|\ge1/2$. Thus
\[
 \frac{|G\cap B(y,at)|}{|B(y,at)|}\le1\le2\frac{|E\cap B(y,bt)|}{|B(y,bt)|}.
\]
Tonelli's theorem gives $\int_G(S_aF)^2\le2\int_E(S_bF)^2$. Using $\Ph(F)\le\int_GF^2+\lambda^2|G^c|$ and $B=\{S_bF>\lambda\}$ proves \eqref{eq:aperture}.
\end{proof}
\begin{lemma}\label{le:osc}
Let $v$ be a Hilbert-valued harmonic function on $\R_+^{d+1}$. If $y\in B(x,at)$ and $\eta>a$, then
\begin{equation}
 \|v(y,t)-v(x,t)\|_\HH\le CA_\eta v(x),\label{eq:osc}
\end{equation}
where $A_\eta v$ is the one-parameter area integral function.
\end{lemma}
\begin{proof}
Choose $c>0$, depending only on $a,\eta$, sufficiently small that every ball of radius $ct$ centered on the horizontal segment from $(x,t)$ to $(y,t)$ lies in $\Gamma_\eta(x)$. By subharmonicity of $\|\nabla v\|_\HH^2$ and the mean-value inequality on balls,
\[
 \|\nabla v(z,t)\|_\HH^2\le Ct^{-d-1}\int_{B((z,t),ct)}\|\nabla v(w,s)\|_\HH^2\,dw\,ds
 \le Ct^{-2}(A_\eta v(x))^2.
\]
Integrating along the horizontal segment proves the assertion.
\end{proof}
We now give the one-parameter Hilbert-valued results required in the proofs of Theorems~\ref{th:local} and~\ref{th:main}. For $q\in L^2(\R^d;\HH)$, write $v(x,t)=P_tq(x)$ and define
\begin{align*}
 &N^+q(x)=\sup_{t>0}\|v(x,t)\|_\HH,\\
 &N_aq(x)=\sup_{t>0,\ |y-x|<at}\|v(y,t)\|_\HH,\\
 &A_bq(x)=\left(\int_{\Gamma_b(x)}\|\nabla v(y,t)\|_\HH^2
       \frac{t\,dy\,dt}{v_d(bt)^d}\right)^{1/2}.
\end{align*}
\begin{lemma}\label{le:local1}
For fixed $b>0$, there exists $\delta>0$ such that, for any measurable $E\subset\R^d$ with $|E^c|<\infty$, the corresponding set $G_E=\{x:M\ind_{E^c}(x)\le\delta\}$ satisfies
\begin{equation}
 |G_E^c|\le C|E^c|,\qquad \int_{G_E}(N^+q)^2\,dx\le C\int_E(A_bq)^2\,dx,\label{eq:local1}
\end{equation}
with constants independent of the Hilbert space $\HH$ and $b, q, E$.
\end{lemma}
\begin{proof}
Put $W_b(y,t)=B(y,bt)$ and $\rho_E(y,t)=|E\cap W_b(y,t)|/|W_b(y,t)|$. Tonelli's theorem yields
\[
 \int_E(A_bq)^2\,dx=\int_{\R^d}\int_0^\infty\|\nabla v(y,t)\|_\HH^2\rho_E(y,t)t\,dt\,dy.
\]
Since $\|\partial_tv\|_\HH\le\|\nabla v\|_\HH$, setting $H=\ind_{\{\rho_E\ge1/2\}}\partial_tv$ gives $\|H\|_{\Ed(\HH)}^2\le2\int_E(A_bq)^2\,dx$. For $G_E=\{M\ind_{E^c}\le\delta\}$, the weak type $(1,1)$ bound of $M$ gives $|G_E^c|\le C\delta^{-1}|E^c|$. If $|y-x|<t/2$, then $B(y,bt)\subset B(x,(b+1/2)t)$, so
\[
 1-\rho_E(y,t)\le C_bM\ind_{E^c}(x).
\]
Taking $\delta\le(2C_b)^{-1}$ ensures $\rho_E(y,t)\ge1/2$. The support of the reconstruction kernel $\psi_{t-s}$ satisfies $|y-x|<(t-s)/2<t/2$. Hence, for $x\in G_E$, all derivative evaluation points in the reconstruction lie in $\{\rho_E\ge1/2\}$. By the one-parameter reconstruction formula (Lemma~\ref{le:reconstruct}) and Lemma~\ref{le:max},
\[
 \int_{G_E}(N^+q)^2\,dx
 \le C\int_{G_E}\sup_{s>0}\|T_sH(x)\|_\HH^2\,dx
 \le C\|H\|_{\Ed(\HH)}^2.
\]
This proves the lemma.
\end{proof}
In Lemma~\ref{le:local1}, take $E=\{A_bq\le\lambda\}$. The definition of $\Ph$ gives
\begin{equation}
 \Ph(N^+q)\le C\left[\int_{\{A_bq\le\lambda\}}(A_bq)^2\,dx
       +\lambda^2|\{A_bq>\lambda\}|\right]=C\Ph(A_bq).\label{eq:radial1}
\end{equation}
\begin{proposition}\label{pr:main1}
For every $q\in L^2(\R^d;\HH)$ and $a,b>0$,
\begin{equation}
 \Ph(N_aq)\le C\Ph(A_bq),\label{eq:main1}
\end{equation}
where $C$ is independent of $\HH$, $\lambda$, and $q$, and depends only on $a,b$ and the dimension $d$.
\end{proposition}
\begin{proof}
Take $\eta>a$. Lemma~\ref{le:osc} gives $N_aq\le N^+q+CA_\eta q$. Applying Lemma~\ref{le:trunc}, \eqref{eq:radial1}, and Lemma~\ref{le:aperture} in succession yields
\[
 \Ph(N_aq)\le C\{\Ph(N^+q)+\Ph(A_\eta q)\}
 \le C\Ph(A_\eta q)\le C\Ph(A_bq).
\]
This proves the proposition.
\end{proof}
The preceding proposition yields the one-parameter version of \eqref{eq:reverse}. Thus Proposition~\ref{pr:main1} gives a new proof of that inequality, avoiding Green's formula and the treatment of sawtooth regions. We next extend this approach to the multiparameter setting, treating the two-parameter case as an example.

\section{Proofs of the multiparameter theorems}
Building on the one-parameter Hilbert-valued estimates in Section~2, we establish a product operator estimate, a product reconstruction formula, and geometric lemmas for the product setting, and then prove the main theorems.

\subsection{Two-parameter reconstruction and product geometric lemmas}
Apply the preceding constructions in each spatial factor. Write $\psi_{i,r}=\psi_{d_i,r}$ and $\chi_{i,r}=\chi_{d_i,r}$, with corresponding operators $S_i,T_{i,s},E_{i,s}$, positive kernels $K_{i,t}$, and averaging operators $V_{i,s}h=\chi_{i,s}*_ih$. For a Hilbert space $\HH$, define
\begin{equation}
 \|H\|_{\Enm(\HH)}^2=\int_{(0,\infty)^2}\int_{\R^n\times\R^m}
 t_1t_2|H(y_1,t_1;y_2,t_2)|^2\,dy_1\,dy_2\,dt_1\,dt_2.\label{eq:energy2}
\end{equation}
When no confusion can arise, we abbreviate $\|H\|_\HH$ as $|H|$. Here $\mathbb Q_+$ denotes the positive rational numbers.
\begin{proposition}\label{pr:max2}
For every Hilbert space $\HH$ and $H\in\Enm(\HH)$,
\begin{equation}
 \left\|\sup_{s_1>0,s_2>0}|T_{1,s_1}T_{2,s_2}H(\cdot)|\right\|_{L^2(X)}^2
 \le C\|H\|_{\Enm(\HH)}^2.\label{eq:max2}
\end{equation}
\end{proposition}
\begin{proof}
By density, we may assume that $H$ is smooth and compactly supported in $\R^n\times\R^m\times(0,\infty)^2$. Operators acting in different variables commute, and
\begin{align*}
 T_{1,s_1}T_{2,s_2}
 &=V_{1,s_1}V_{2,s_2}S_1S_2+E_{1,s_1}V_{2,s_2}S_2\\
 &\quad+V_{1,s_1}S_1E_{2,s_2}+E_{1,s_1}E_{2,s_2}.
\end{align*}
Denote the pointwise suprema of the $\HH$-norms of these four terms by $U_{11},U_{10},U_{01},U_{00}$, respectively. The first term satisfies
\[
 U_{11}\le CM_1M_2|S_1S_2H|,\qquad
 \|U_{11}\|_2^2\le C\|S_1S_2H\|_2^2\le C\|H\|_{\Enm(\HH)}^2.
\]
In the second term, $E_{1,s_1}$ and $V_{2,s_2}$ act in different variables and hence commute, so
\[
 U_{10}\le CM_2\left(\sup_{s_1>0}|E_{1,s_1}S_2H|\right).
\]
By \eqref{eq:Ebound}, Fubini's theorem, and \eqref{eq:Sbound} applied in the second variable for each fixed $(y_1,t_1)$,
\begin{align*}
 \left\|\sup_{s_1>0}|E_{1,s_1}S_2H|\right\|_2^2
 &\le C\int_{\R^n}\int_0^\infty
 t_1\|(S_2H)_{t_1}(y_1,\cdot)\|_{L^2(\R^m;\HH)}^2\,dt_1\,dy_1\\
 &\le C\|H\|_{\Enm(\HH)}^2.
\end{align*}
The $L^2$ boundedness of $M_2$ therefore gives $\|U_{10}\|_2^2\le C\|H\|_{\Enm(\HH)}^2$. The third term is treated symmetrically.

For the final term, apply \eqref{eq:Ebound} successively, first in the first variable and then in the second, moving the suprema over the time variables inside the nonnegative integrals. We obtain
\begin{align*}
 U_{00}(x_1,x_2)^2
 &\le C\int_0^\infty\int_0^\infty t_1t_2
 \int_{\R^n}\int_{\R^m}K_{1,t_1}(x_1-y_1)K_{2,t_2}(x_2-y_2)\\
 &\hspace{35mm}\times|H(y_1,t_1;y_2,t_2)|^2\,dy_2\,dy_1\,dt_2\,dt_1.
\end{align*}
Integrating in $x_1,x_2$ and using $\|K_{i,t}\|_1=\|k_i\|_1$ yields the $L^2$ estimate for the double error term. Combining all four estimates proves \eqref{eq:max2}.
\end{proof}
We next give the two-parameter reconstruction formula for functions harmonic in each group of variables.
\begin{lemma}\label{le:reconstruct2}
For the product Poisson extension $u(y_1,t_1;y_2,t_2)=P_{t_1}^{(1)}P_{t_2}^{(2)}f(y_1,y_2)$, we have
\begin{equation}
 u(y_1,s_1;y_2,s_2)=(c_nc_m)^{-1}T_{1,s_1}T_{2,s_2}
 (\partial_{t_1}\partial_{t_2}u(\cdot,t_1;\cdot,t_2))(y_1,y_2).\label{eq:reconstruct2}
\end{equation}
The formula on the right uses only values of $u$ at points satisfying $t_i>s_i$ and $|y_i-x_i|<(t_i-s_i)/2<t_i/2$.
\end{lemma}
\begin{proof}
First apply Lemma~\ref{le:reconstruct} in the first spatial variable, treating the second spatial variable and its time variable as Hilbert-valued parameters. Apply the same lemma in the second spatial variable. Multiplying the constants from the two one-parameter reconstructions and canceling the two minus signs gives \eqref{eq:reconstruct2}. Since $\psi_{i,t_i-s_i}$ is supported in $B(0,(t_i-s_i)/2)$, the reconstruction integral uses only values of $u$ in the stated cones.
\end{proof}
\begin{lemma}\label{le:densityid}
For every measurable set $E\subset X$,
\begin{equation}
 \int_E(A_\beta f)^2\,dx=\int_{(0,\infty)^2}\int_X
 |\grad_1\grad_2u(y,t)|^2\rho_E(y,t)t_1t_2\,dy\,dt,\label{eq:densityid}
\end{equation}
where
\[
 \rho_E(y,t)=\rho_E^\beta(y,t)=\frac{|E\cap W_\beta(y,t)|}{|W_\beta(y,t)|},\quad
 W_\beta(y,t)=B(y_1,\beta_1t_1)\times B(y_2,\beta_2t_2),
\]
and $\beta\in(0,\infty)^2$. The identity remains valid if $\grad_1\grad_2u(y,t)$ is replaced by an arbitrary measurable function and $A_\beta$ by the square function $S_\beta$ in Lemma~\ref{le:aperture2} below.
\end{lemma}
\begin{proof}
Substitute the definition of the area integral function and apply Tonelli's theorem. The details are omitted.
\end{proof}
\begin{lemma}\label{le:geometry}
Suppose $a_i,\beta_i>0$. There exists $\delta>0$, depending only on these numbers and the dimensions, such that for every $E\subset X$ with $|E^c|<\infty$, the set $G_E=\{x:M_1M_2\ind_{E^c}(x)\le\delta\}$ satisfies
\begin{equation}
 \begin{gathered}
 |G_E^c|\le C|E^c|,\\
 \rho_E(y,t)\ge\tfrac12\quad\text{if }x\in G_E,\quad |y_i-x_i|<a_it_i\ (i=1,2).
 \end{gathered}\label{eq:geometry}
\end{equation}
\end{lemma}
\begin{proof}
Let $B=E^c$. By the $L^2$ boundedness of $M_1M_2$,
\[
 |\{M_1M_2\ind_B>\delta\}|\le\delta^{-2}\|M_1M_2\ind_B\|_2^2\le C\delta^{-2}|B|.
\]
For the second estimate, put $B_i^*=B(x_i,(\beta_i+a_i)t_i)$. Then $B(y_i,\beta_it_i)\subset B_i^*$, so
\[
 1-\rho_E(y,t)\le\frac{C_{\beta,a}}{|B_1^*||B_2^*|}
 \int_{B_1^*\times B_2^*}\ind_B(w)\,dw
 \le C_{\beta,a}M_1M_2\ind_B(x).
\]
Taking $\delta\le(2C_{\beta,a})^{-1}$ proves the second estimate.
\end{proof}
\begin{lemma}\label{le:aperture2}
Let $F$ be a Hilbert-valued measurable function on the product space. Define
\[
 S_cF(x)=\left(\iint_{\Gamma_c(x)}\|F(y,t)\|_\HH^2
 \prod_{i=1}^2\frac{t_i\,dy_i\,dt_i}{v_{d_i}(c_it_i)^{d_i}}\right)^{1/2}.
\]
Then
\begin{equation}
 \Ph(S_aF)\le C\Ph(S_bF),\label{eq:aperture2}
\end{equation}
where $C$ depends only on $a,b$ and the dimensions.
\end{lemma}
\begin{proof}
The assertion is immediate if the right-hand side is infinite. Otherwise, set $E=\{S_bF\le\lambda\}$, $B=E^c$, and $G=\{M_1M_2\ind_B\le\varepsilon\}$. Then $|G^c|\le C_\varepsilon|B|$. If $G\cap W_a(z)$ is nonempty, choose a point in this intersection as a base point. By Lemma~\ref{le:geometry}, for a sufficiently small fixed $\varepsilon$ we have $\rho_E^b(z)\ge1/2$, and hence
\[
 \frac{|G\cap W_a(z)|}{|W_a(z)|}\le2\frac{|E\cap W_b(z)|}{|W_b(z)|}.
\]
The inequality also holds when its left-hand side is zero. Tonelli's theorem gives $\int_G(S_aF)^2\le2\int_E(S_bF)^2$. Splitting the integral defining $\Ph$ into the parts over $G$ and $G^c$ proves \eqref{eq:aperture2}.
\end{proof}

\subsection{Proof of Theorem \ref{th:local}}
\begin{proof}[Proof of Theorem~\ref{th:local}]
By a limiting argument, we may assume that $f\in C_c^\infty(\R^n\times\R^m)$ and $u=P_{t_1}^{(1)}P_{t_2}^{(2)}f$. In Lemma~\ref{le:geometry}, take $a_1=a_2=1/2$ and choose $G_E$ accordingly. Then
\[
 \rho_E(y,t)\ge\tfrac12\quad\text{whenever }x\in G_E,\quad |y_i-x_i|<t_i/2.
\]
Set
\[
 \Omega_E=\{(y,t):\rho_E(y,t)\ge1/2\},\qquad
 H=\ind_{\Omega_E}\partial_{t_1}\partial_{t_2}u.
\]
Since $\partial_{t_1}\partial_{t_2}$ is one component of $\grad_1\grad_2$, Lemma~\ref{le:densityid} gives
\begin{equation}
 \|H\|_{\Enm}^2\le2\int_E(A_\beta f)^2\,dx.\label{eq:Hlocal}
\end{equation}
For $x\in G_E$, Lemma~\ref{le:reconstruct2} and the preceding geometric property show that every point used in the reconstruction belongs to $\Omega_E$. Thus $\partial_{t_1}\partial_{t_2}u$ can be replaced by $H$ in the reconstruction integral. Proposition~\ref{pr:max2} gives
\begin{align*}
 \int_{G_E}(N^+f)^2\,dx
 &\le C\int_{G_E}\sup_{s_1>0,s_2>0}|T_{1,s_1}T_{2,s_2}H|^2\,dx\\
 &\le C\|H\|_{\Enm}^2\le C\int_E(A_\beta f)^2\,dx.
\end{align*}
This proves Theorem~\ref{th:local}.
\end{proof}
Lemma~\ref{le:local1} gives the corresponding one-parameter local estimate for Hilbert-valued functions.
\begin{corollary}\label{co:radial2}
Under the hypotheses of Theorem~\ref{th:local}, for every $f\in L^2(\R^n\times\R^m)$ and $\lambda>0$,
\begin{equation}
 \Ph(N^+f)\le C\Ph(A_\beta f).\label{eq:radial2}
\end{equation}
\end{corollary}
\begin{proof}
Take $E=\{A_\beta f\le\lambda\}$, so that $|E^c|<\infty$. Split the integral defining $\Ph(N^+f)$ over $G_E$ and $G_E^c$, and apply Theorem~\ref{th:local}.
\end{proof}

\subsection{Proof of Theorem \ref{th:main}}
\begin{lemma}\label{le:Q}
For $\eta_i>\alpha_i$, $f\in C_c^\infty(\R^n\times\R^m)$, $u=P_{t_1}^{(1)}P_{t_2}^{(2)}f$, and $z_i=(y_i,t_i)$, define
\begin{align*}
 Q_1(x)&=\sup_{z_1\in\Gamma_{1,\alpha_1}(x_1)}
 \left(\int_{\Gamma_{2,\eta_2}(x_2)}|\grad_2u(z_1,z_2)|^2
 \frac{t_2\,dy_2\,dt_2}{v_m(\eta_2t_2)^m}\right)^{1/2},\\
 Q_2(x)&=\sup_{z_2\in\Gamma_{2,\alpha_2}(x_2)}
 \left(\int_{\Gamma_{1,\eta_1}(x_1)}|\grad_1u(z_1,z_2)|^2
 \frac{t_1\,dy_1\,dt_1}{v_n(\eta_1t_1)^n}\right)^{1/2}.
\end{align*}
Then
\begin{equation}
 \Ph(Q_1)+\Ph(Q_2)\le C\Ph(A_\beta f).\label{eq:Q}
\end{equation}
\end{lemma}
\begin{proof}
We first estimate $Q_1$. Fix $\eta_2$, and let
\[
 d\mu_{2,\eta_2}(y_2,t_2)=\frac{t_2\,dy_2\,dt_2}{v_m(\eta_2t_2)^m},\qquad
 \HH_2=L^2\bigl(\R^m\times(0,\infty),d\mu_{2,\eta_2};\C^{m+1}\bigr).
\]
This is a fixed Hilbert space independent of $x_2$. For each $x_2$, define the $\HH_2$-valued function
\[
 q_{x_2}(y_1)(y_2,t_2)=\ind_{\{|y_2-x_2|<\eta_2t_2\}}\grad_2P_{t_2}^{(2)}f(y_1,y_2).
\]
Fubini's theorem and the Fourier transform give
\begin{equation}
 \begin{split}
 \int_{\R^m}\|q_{x_2}\|_{L^2(\R^n;\HH_2)}^2\,dx_2
 &=\int_{\R^n}\int_{\R^m\times(0,\infty)}
 t_2|\grad_2P_{t_2}^{(2)}f(y_1,y_2)|^2\,dy_2\,dt_2\,dy_1\\
 &=\tfrac12\|f\|_2^2.
 \end{split}\label{eq:qnorm}
\end{equation}
Thus $q_{x_2}\in L^2(\R^n;\HH_2)$ for almost every $x_2$. The Hilbert-valued nontangential maximal function of its Poisson extension in the first variable is $Q_1(\cdot,x_2)$, and the corresponding area integral function is $A_{(\beta_1,\eta_2)}f(\cdot,x_2)$. Applying Proposition~\ref{pr:main1} for each $x_2$, integrating, and then using Lemma~\ref{le:aperture2}, we obtain
\[
 \Ph(Q_1)\le C\Ph(A_{(\beta_1,\eta_2)}f)\le C\Ph(A_\beta f).
\]
The estimate for $Q_2$ is symmetric: fix $x_1$ and use the one-parameter argument with values in
\[
 \HH_1=L^2\left(\R^n\times(0,\infty),
 \frac{t_1\,dy_1\,dt_1}{v_n(\eta_1t_1)^n};\C^{n+1}\right).
\]
This gives $\Ph(Q_2)\le C\Ph(A_\beta f)$.
\end{proof}
\begin{proof}[Proof of Theorem~\ref{th:main}]
We may assume that $f\in C_c^\infty(\R^n\times\R^m)$ and $u$ is its two-parameter Poisson extension. Consider the decomposition
\begin{align*}
 u(z_1,z_2)-u(x_1,t_1;x_2,t_2)
 &=\bigl(u(z_1,z_2)-u(z_1;x_2,t_2)\bigr)\\
 &\quad+\bigl(u(z_1;x_2,t_2)-u(x_1,t_1;x_2,t_2)\bigr).
\end{align*}
For the first difference, fix $z_1$, regard $u(z_1,z_2)$ as a harmonic function in the second group of variables, and apply Lemma~\ref{le:osc} with $a=\alpha_2$ and $\eta=\eta_2$. Its right-hand side is controlled by $Q_1(x)$. For the second difference, fix $z_2=(x_2,t_2)$ and apply the same lemma in the first group of variables with $a=\alpha_1$ and $\eta=\eta_1$; its right-hand side is controlled by $Q_2(x)$. Taking suprema over $z_1,z_2$ gives
\begin{equation}
 N_\alpha f\le N^+f+C(Q_1+Q_2).\label{eq:decomp}
\end{equation}
By Lemma~\ref{le:trunc}, \eqref{eq:radial2}, and Lemma~\ref{le:Q},
\[
 \Ph(N_\alpha f)\le C\{\Ph(N^+f)+\Ph(Q_1)+\Ph(Q_2)\}\le C\Ph(A_\beta f).
\]
This proves Theorem~\ref{th:main}.
\end{proof}
\begin{corollary}\label{co:Lp}
Under the hypotheses of Theorem~\ref{th:main},
\begin{equation}
 |\{N_\alpha f>\lambda\}|\le C|\{A_\beta f>\lambda\}|
 +C\lambda^{-2}\int_{\{A_\beta f\le\lambda\}}(A_\beta f)^2\,dx.\label{eq:distribution}
\end{equation}
Moreover, for every $0<p<2$,
\begin{equation}
 \|N_\alpha f\|_{L^p}\le C_p\|A_\beta f\|_{L^p},\qquad
 \|N_\alpha f\|_{L^{p,\infty}}\le C_p\|A_\beta f\|_{L^{p,\infty}}.\label{eq:Lp}
\end{equation}
\end{corollary}
\begin{proof}
Since $\lambda^2|\{N_\alpha f>\lambda\}|\le\Ph(N_\alpha f)\le C\Ph(A_\beta f)$, splitting the right-hand side over $\{A_\beta f\le\lambda\}$ and $\{A_\beta f>\lambda\}$ gives \eqref{eq:distribution}. The strong type assertion follows from Tonelli's theorem and
\[
 \int_0^\infty\lambda^{p-3}\min\{s^2,\lambda^2\}\,d\lambda
 =\frac{2s^p}{p(2-p)},\qquad0<p<2.
\]
For the weak type assertion, writing $K=\|F\|_{L^{p,\infty}}$, we have
\[
 \Ph(F)=\int_0^\lambda2r|\{F>r\}|\,dr\le\frac{2}{2-p}K^p\lambda^{2-p}.
\]
Take $F=A_\beta f$ and apply Theorem~\ref{th:main}.
\end{proof}
\begin{remark}
For any fixed finite number of factors, one can formally expand the tensor product of the operators $V_iS_i+E_i$. Each term must be treated separately using the corresponding positive integrable kernel estimates and iterated maximal operator estimates. Thus estimate \eqref{eq:reverse} extends to the multiparameter setting.
\end{remark}

\par\medskip
\noindent\textbf{Acknowledgments.} The author was supported by the National Natural Science Foundation of China (No.~12571109) and the National Key Research and Development Program of the Ministry of Science and Technology of China (No.~2022YFA1005700).

\end{document}